\documentclass[12pt,a4paper]{article}
\usepackage{amssymb}

\newtheorem{defn}{Definition}
\newtheorem{thm}{Theorem}
\newtheorem{col}{Corollary}
\newtheorem{lem}{Lemma}

\newenvironment{proof}{\begin{trivlist} \item[] {\em Proof}. }%
{\hfill $\square$ \end{trivlist}}

\def\sn{\mathop{\rm sn}\nolimits}
\def\cn{\mathop{\rm cn}\nolimits}
\def\dn{\mathop{\rm dn}\nolimits}

\title{
THE GENERAL ANALYTIC SOLUTION OF A FUNCTIONAL EQUATION OF ADDITION TYPE\\}

\author{H. W. Braden\\
\normalsize
\em Department of Mathematics and Statistics,\\
\normalsize
\em The University of Edinburgh, \\
\normalsize
\em Edinburgh, UK. \\
\normalsize
e-mail: hwb@ed.ac.uk
\\
V. M. Buchstaber\\
\normalsize
   \em      National Scientific and Research Institute of\\
\normalsize
     \em    Physico-Technical and Radio-Technical Measurement,\\
\normalsize
    \em     VNIIFTRI, Mendeleevo, Moscow Region\\
\normalsize
      \em   141570, Russia\\
\normalsize
e-mail: buchstab@tomogr.msk.su
\\
\normalsize
{\tt funct-an/9508002}\\
}

\begin{document}

\renewcommand{\thepage}{}
\begin{titlepage}
\maketitle

\begin{abstract}
The general analytic solution to the functional equation
$$
\phi_1(x+y)=
{ { \biggl|\matrix{\phi_2(x)&\phi_2(y)\cr\phi_3(x)&\phi_3(y)\cr}\biggr|}
\over
{ \biggl|\matrix{\phi_4(x)&\phi_4(y)\cr\phi_5(x)&\phi_5(y)\cr}\biggr|} }
$$
is characterised. Up to the action of the symmetry group, this is described
in terms of Weierstrass elliptic functions.
We illustrate our theory  by applying it to the classical addition theorems
of the Jacobi elliptic functions and the  functional equations
$$
\phi_1(x+y)=\phi_4(x)\phi_5(y)+\phi_4(y)\phi_5(x)
$$
and
\[
\Psi _1(x+y)=\Psi _2(x+y) \phi_2(x)\phi_3(y) +\Psi_3(x+y) \phi_4(x)\phi_5(y).
\]
\end{abstract}
\vfill
\end{titlepage}
\renewcommand{\thepage}{\arabic{page}}

\section{Introduction}
The purpose of this article is to describe the general analytic solution to
the functional equation
\begin{equation}
\phi_1(x+y)=
{ { \biggl|\matrix{\phi_2(x)&\phi_2(y)\cr\phi_3(x)&\phi_3(y)\cr}\biggr|}
\over
{ \biggl|\matrix{\phi_4(x)&\phi_4(y)\cr\phi_5(x)&\phi_5(y)\cr}\biggr|} }.
\label{functional}
\end{equation}
Although this equation appears to depend on five a priori unknown functions
we shall show that (\ref{functional}) is invariant under a 
large group of symmetries ${\cal{G}}$ and that each orbit has a 
solution of a particularly nice form, expressible in terms elliptic functions:
\begin{thm}
The general analytic solution to the
functional equation (\ref{functional})\ is, up to a ${\cal{G}}$ action given
by (\ref{symm}-\ref{symmsb}), of the form
$$\phi_1(x)= {\Phi(x;\nu_1)\over\Phi(x;\nu_2)},\quad
  {\phi_2(x)\choose\phi_3(x)}={\Phi(x;\nu_1)\choose\Phi(x;\nu_1)\sp\prime}\quad
{ and}\quad
  {\phi_4(x)\choose\phi_5(x)}={\Phi(x;\nu_2)\choose\Phi(x;\nu_2)\sp\prime}.
$$
Here
\begin{equation}
\Phi(x;\nu)\equiv {\sigma(\nu-x)\over {\sigma(\nu)\sigma(x)}}\,
  e\sp{\zeta(\nu)x}
\label{soln}
\end{equation}
where $\sigma(x)=\sigma(x|\omega,\omega\sp\prime)$ and
$\zeta(x) ={\sigma(x)\sp\prime \over\sigma(x)}$
are the Weierstrass sigma and zeta functions.
\end{thm}
The group ${\cal G}$ of symmetries of (\ref{functional})\ will
be described  further below. Our proof is constructive and indeed yields
more: 
\begin{thm}
Let $x_0$ be a generic point for (\ref{functional}).
Then (for $k=1,2$) we have
\begin{eqnarray*}
\partial_y \ln
         \left|
          \matrix{\phi_{2k}(x+x_0)&\phi_{2k}(y+x_0)\cr
          \phi_{2k+1}(x+x_0)&\phi_{2k+1}(y+x_0)\cr}
         \right|
\Biggl|_{y=0}
&=&
\zeta(\nu_k)-\zeta(x)-\zeta(\nu_k-x)-\lambda_k
\\
&=&
 -{1\over x}-\lambda_k +\sum_{l=0}F_l\, {x\sp{l+1}\over (l+1)!}
\end{eqnarray*}
and the Laurent expansion
determines the parameters $g_1$, $g_2$ (which are the same for both $k=1,2$)
characterising the elliptic functions of (\ref{soln}) by
$$
g_2={5\over3}\left({F_2+6 F_0\sp2 }\right),
\quad\quad
g_3= 6 F_0\sp3 -F_1\sp2 +{5\over3}F_0 F_2,
$$
and the parameters $\nu_k$ via $F_0=-\wp(\nu_k)$. Further, we have
$$
\phi_1(x+2 x_0)= \phi_1(2 x_0)\,
e\sp{(\lambda_2-\lambda_1) x}\,{\Phi(x;\nu_1)\over\Phi(x;\nu_2)}\,
$$
and
$$
\left( \matrix{\phi_{2k}(x+x_0)\cr \phi_{2k+1}(x+x_0)\cr} \right)
=
{e\sp{-\lambda_k x}\over f(x)}
\Biggl(
\matrix{\phi_{2k}\sp\prime(x_0)&\phi_{2k}(x_0)\cr
        \phi_{2k+1}\sp\prime(x_0)&\phi_{2k+1}(x_0)\cr}
\Biggr)
\Biggl( \matrix{1&0\cr \lambda_k&-1\cr}\Biggr)
\Biggl( \matrix{\Phi(x;\nu_k)\cr \Phi\sp\prime(x;\nu_k)\cr}\Biggr).
$$
Here the function 
$$
f(x)=
{e\sp{-\lambda_k x}\over \Phi(x;\nu_k)}\,
{
\left|
\matrix{\phi_{2k}\sp\prime(x_0)&\phi_{2k}(x_0)\cr
        \phi_{2k+1}\sp\prime(x_0)&\phi_{2k+1}(x_0)\cr}
\right|
\over
\left|
\matrix{\phi_{2k}(x+x_0)&\phi_{2k}(x_0)\cr
        \phi_{2k+1}(x+x_0)&\phi_{2k+1}(x_0)\cr}
\right|
}
$$
is in fact the same for $k=1,2$.
\end{thm}
The term \lq generic\rq\ will be defined below and we will give
more expressions for the quantities appearing in the theorem.

One merit of writing (\ref{functional}) in this general form is that
several different functional equations may now be seen as 
different points on a ${\cal G}$-orbit of
(\ref{functional}). Thus for example
\begin{eqnarray}
\phi_1(x+y)&=&\phi_1(x)\phi_1(y)\label{exponential}\\
\phi_1(x+y)&=&\phi_4(x)\phi_5(y)+\phi_4(y)\phi_5(x)\label{dexponential}\\
A(x+y)[ B(x)-B(y)]&=&A(x) A\sp\prime(y) -A(y)A\sp\prime(x).
\label{calfun}
\end{eqnarray}
are particular\footnote{
These correspond to
\begin{itemize}
\item[(a)]$\phi_2(x)=\phi_1(x)\phi_4(x)$ and $\phi_3(x)=\phi_1(x)\phi_5(x)$,
\item[(b)]$\phi_2(x)=\phi_4\sp2(x)$ and $\phi_3(x)=\phi_5\sp2(x)$,
\item[(c)]$\phi_1(x)=\phi_2(x)=A(x)$, $\phi_3(x)=A\sp\prime(x)$, 
$\phi_4(x)=B(x)$ and $\phi_5(x)=1$.
\end{itemize}
} examples of (\ref{functional}).
The functional equation for the exponential (\ref{exponential})
corresponds to $\nu_1=\nu_2$ in our solution and the exponential comes
wholly from ${\cal G}$.
Particular cases of (\ref{dexponential}) have been studied in
\cite{BK} and we shall determine (see Lemma 4)
the general solution to (\ref{dexponential}) as an application of our work.

More interesting is  equation (\ref{calfun}) which has been
studied by various authors
with assumptions of even/oddness on the functions appearing
\cite{Ca2, OPc, PS} or assumptions on the nature of $B$ \cite{Kr1}.
The general solution \cite{BCb, Bu1} $A(x)=\Phi(x;\nu)$ now
corresponds to to the limit $\nu_2\rightarrow0$ together with a
${\cal G}$ action. This will be illustrated later.  

Finally, when
$\phi_1(x)=\alpha(x)$, $\phi_2(x)=\alpha(x)\tau(x)$, $\phi_4(x)=\tau(x)$,
$\phi_3(x)=\phi_2\sp\prime(x)$ and $\phi_5(x)=\phi_4\sp\prime(x)$
we obtain the functional equation
\begin{equation}
\alpha(x+y)=\alpha(x)\alpha(y)+\tau(x)\tau(y)\psi(x+y).
\label{calogero}
\end{equation}
The function $\psi(x)$ will be described in more detail in the sequel.
This equation was studied by Bruschi and Calogero \cite{BCb}
and will be used in our analysis.

Let us remark that both (\ref{dexponential})  and (\ref{calogero})
may be viewed as limiting cases of the functional equation
\begin{equation}
\Psi _1(x+y)=\Psi _2(x+y) \phi_2(x)\phi_3(y) +\Psi_3(x+y) \phi_4(x)\phi_5(y),
\label{eq:biggy}
\end{equation}
which \textit{a priori } depends on seven unknown
functions. Later we shall show how (\ref{functional}) may be used to solve
this.

It remains to place (\ref{functional}) in some form of context.
The last decade has seen a remarkable confluence of ideas from
completely integrable systems, geometry, field theory and functional
equations that is still being assimilated. To make  some of these matters
concrete let us consider how such functional equations arise in the
context of integrable systems of particles on the line.
A pair of matrices ${L,M}$ such that ${\dot L}= [L,M]$ is known as a
Lax Pair; this is a zero curvature condition. 
Starting with an ansatz for the matrices $L$ and $M$ one
seeks  restrictions  necessary to obtain equations of motion of  some desired
form. These restrictions typically involve the study of functional equations.
The paradigm for this approach is the Calogero-Moser \cite{Ca}
system.
Beginning with the ansatz (for $n\times n$ matrices)
$$L_{jk}=p_j\delta_{jk}+ \, g\,(1-\delta_{jk}) A(q_j-q_k),\quad
  M_{jk}=g\, [\delta_{jk}\sum_{l\ne j}B(q_j-q_l)-(1-\delta_{jk})C(q_j-q_k) ]
$$
%( with $b(x)$, $c(x)$ even functions and $a(x)$ an odd function)
one finds ${\dot L}= [L,M]$
yields the
equations of motion for the Hamiltonian system ($n\ge3$)
$$H={1\over 2}\sum_{j}p_j\sp2 +g\sp2\sum_{j<k}U(q_j-q_k)
\quad\quad\quad U(x)=A(x) A(-x) + {\rm constant}$$
provided $C(x)=-A\sp\prime(x)$,
and that $A(x)$ and $B(x)$ satisfy the functional equation (\ref{calfun}).
With this ansatz and assuming $B(x)$ even\footnote{
This assumption can in fact be removed \cite{BB}. }
Calogero \cite{Ca2} found
$A(x)$  to be given by (\ref{soln}). In this case the corresponding potential 
is the Weierstrass $\wp$-function: $A(x)A(-x)=\wp(\nu)-\wp(x)$.
The functional equation (\ref{calogero}) is associated with a different
ansatz and yields the relativistic Calogero-Moser systems \cite{BCa, Ra}.
Similarly (\ref{functional}) arises from a more general ansatz \cite{BB}
associated with equations of motion of the form
$$
{\ddot q_j}=\sum_{k\ne j}(a+b \dot q_j) (a+b \dot q_k) 
V_{jk}(q_j-q_k),
$$
which combines both relativistic ($b\ne0$) and nonrelativistic
($b=0$) systems together with potentials that can vary between particle pairs.
This unifies, for example, Calogero-Moser and Toda systems \cite{RS, Ra, Rb}. 
The relativistic
examples yield the functional equation (\ref{functional}) while the
nonrelativistic situation involves the functional equation
\begin{equation}
\phi_6(x+y)=\phi_1(x+y)\big( \phi_4(x)-\phi_5(y)\big)
+{ \biggl|\matrix{\phi_2(x)&\phi_3(y)\cr\phi\sp\prime_2(x)&\phi\sp\prime_3(y)
\cr}\biggr|}.
\label{functional2}
\end{equation}
The general analytic solution to (\ref{functional2}) has yet to
be determined although particular solutions are known.
We remark that (\ref{functional}), and after suitably
symmetrizing (\ref{functional2}), are
particular cases of the functional equation
\begin{equation}
\sum_{i=0}\sp{N}
\phi_{3i}(x+y)
 { \biggl|\matrix{\phi_{3i+1}(x)&\phi_{3i+1}(y)\cr
                  \phi_{3i+2}(x)&\phi_{3i+2}(y)\cr}\biggr|}=0
\label{Bfun}
\end{equation}
with $N=1$ in the former case and $N=2$ in  the latter.
When $\phi_{3i+2}=\phi_{3i+1}\sp\prime$
Buchstaber and Krichever have discussed  (\ref{Bfun}) in connection with
functional equations satisfied by Baker-Akhiezer functions \cite{BKr}.

Lax pairs are but one way in which functional equations are associated
with integrable systems and we mention \cite{BFV, Ca3, BP, GT} for others.
There also appears a close connection between these functional
equations and the elliptic genera associated with the 
string inspired Witten index \cite{HBJ, Kr}. Krichever for example
used the functional equation (\ref{calfun}) in his proof
of the \lq rigidity\rq\ property of elliptic genera \cite{Kr}
and it also  appears when discussing rational and pole solutions
of the KP and KdV equations \cite{Kr1, AMM}.
We feel this connection between functional equations and
completely integrable systems is part of a broader and less well
understood aspect of the subject that deserves further attention.

An outline of the paper is as follows. First we will discuss the group of
symmetries of (\ref{functional}). These will be used in the proof of
theorem 1. Before turning to the proof we show in section 3
how the indicated solution
does indeed satisfy (\ref{functional}), using this as a vehicle to recall
some of the properties of elliptic functions that we will need throughout.
Section 4 is devoted to the proof of theorem 1 and
section 5 to that of theorem 2. Several applications of our theorems
including the general analytic solution to (\ref{dexponential})
and a discussion of (\ref{eq:biggy})
are then given in section 6. An appendix is given containing
various elliptic function formulae we shall make use of.

Various versions of Theorem 1 have appeared in unpublished preprints.
In \cite{Bu2} the form of $\phi_1(x)$ only was stated. In
\cite{BB} we introduced the $\cal G$ action to give Theorem 1 in
its present form. In improving the proof of this we obtained Theorem 2, 
given here alongside the better proof of Theorem 1.

\section{The Group of Symmetries}
We next describe the  group ${\cal G}$ of
invariances  of (\ref{functional}).
Theorem $1$ gives a representative of each  ${\cal G}$  orbit
on the solutions of (\ref{functional})\ with a particularly nice form.
First observe that a large group of symmetries ${\cal G}$ act on the solutions
of (\ref{functional}). The transformation
\begin{equation}
\Biggl( \phi_1(x), {\phi_2(x)\choose\phi_3(x)},{\phi_4(x)\choose\phi_5(x)} 
\Biggr) \rightarrow
  \Biggl( c\, e\sp{\lambda x} \phi_1(x), U {e\sp{-\lambda\sp\prime x}
\phi_2(x)\choose e\sp{-\lambda\sp\prime x}\phi_3(x)},
V {e\sp{\lambda\sp{\prime\prime}x} \phi_4(x)\choose
   e\sp{\lambda\sp{\prime\prime}x} \phi_5(x)} \Biggr)
\label{symm}
\end{equation}
clearly preserves (\ref{functional})\ provided
\begin{equation}
\lambda+\lambda\sp\prime+\lambda\sp{\prime\prime}=0,\quad\quad
  U,V\in GL_2,\quad\quad {\rm and}\quad\quad \det U=c\,\det V.
\label{constraints}
\end{equation}
Further, (\ref{functional}) is also preserved by
\begin{equation}
\Biggl( \phi_1(x), {\phi_2(x)\choose\phi_3(x)},{\phi_4(x)\choose\phi_5(x)}
\Biggr) \rightarrow\Biggl(
{1\over \phi_1(x)},{\phi_4(x)\choose\phi_5(x)},{\phi_2(x)\choose\phi_3(x)}
\Biggr)
\label{symms}
\end{equation}
and
\begin{equation}
\Biggl( \phi_1(x), {\phi_2(x)\choose\phi_3(x)},{\phi_4(x)\choose\phi_5(x)}
\Biggr) \rightarrow\Biggl(\phi_1(x),
f(x){\phi_2(x)\choose\phi_3(x)},f(x){\phi_4(x)\choose\phi_5(x)}\Biggr).
\label{symmsb}
\end{equation}
We will use these symmetries in our proof of theorem $1$ to find a
solution of (\ref{functional})\ on each  ${\cal G}$  orbit
with a particularly nice form.

\section{Illustration of the Solution}
Before proceeding to the proof it is instructive
to see how the stated solution satisfies (\ref{functional}).
This will also allow us to introduce some elliptic function identities needed
throughout.
From the definition of the zeta function we have
\begin{equation}
\bigr(\ln \Phi(x;\nu)\bigl)\sp\prime = -\zeta(\nu-x)-\zeta(x)+\zeta(\nu).
\label{defzeta}
\end{equation}
Thus
\begin{eqnarray*}
 {\biggl| \matrix{\Phi(x;\nu)&\Phi(y;\nu)\cr\Phi(x;\nu)\sp\prime&
     \Phi(y;\nu)\sp\prime\cr}\biggr| }
&=&{\Phi(x;\nu)\Phi(y;\nu)\biggl[ \bigr(\ln \Phi(y;\nu)\bigl)\sp\prime
  -\bigr(\ln \Phi(x;\nu)\bigl)\sp\prime\biggr] } \hfill\\
&=&\Phi(x;\nu)\Phi(y;\nu)\biggl[
\zeta(\nu-x)+\zeta(x)+\zeta(-y)+\zeta(y-\nu)\biggr]. \hfill
\end{eqnarray*}
Upon using the definition of $\Phi$, the right hand side of this
equation takes the form
\begin{equation}
\Phi(x+y;\nu){\sigma(\nu-x)\sigma(\nu-y)\sigma(x+y)\over
    \sigma(\nu-x-y)\sigma(\nu)\sigma(x)\sigma(y) }
    \biggl[ \zeta(\nu-x)+\zeta(x)+\zeta(-y)+\zeta(y-\nu)\biggr].
\label{rhs}
\end{equation}
After noting  the two identities \cite{WW}
\begin{equation}
\zeta(x)+\zeta(y)+\zeta(z)-\zeta(x+y+z)=
{\sigma(x+y)\sigma(y+z)\sigma(z+x)\over\sigma(x)\sigma(y)\sigma(z)\sigma(x+y+z)}
\label{eq:zetas}
\end{equation}
and 
\begin{equation}
\wp (x)-\wp(y)=
{\sigma(y-x)\sigma(y+x)\over\sigma\sp2(y)\sigma\sp2(x)}
\label{eq:wps}
\end{equation}
we find (\ref{rhs})  simplifies to
$\Phi(x+y;\nu)\bigl[{\wp }(x)-{\wp }(y)\bigr] $,
where ${\wp }(x)=-\zeta\sp\prime(x)$ is the Weierstrass ${\wp }$-function. 
Putting these together yields the addition formula
\begin{equation}
\Phi(x+y;\nu)=
{{\biggl| \matrix{\Phi(x;\nu)&\Phi(y;\nu)\cr\Phi(x;\nu)\sp\prime&
     \Phi(y;\nu)\sp\prime\cr}\biggr| }\over
{\wp }(x)-{\wp }(y) },
\label{addn}
\end{equation}
and consequently a solution of (\ref{functional})\ with the stated form. 
Further, from (\ref{addn}) we see the solution to (\ref{calfun})
mentioned in the introduction.

The general solution (\ref{soln}) involves the two nonzero constants
$\nu_1$ and $\nu_2$. 
Let us see how our group of symmetries enables $\phi_1(x)=\Phi(x;\nu_1)$
to occur as a limit $\nu_2\rightarrow0$.
Consider the  ${\cal G}$ action on the general solution $\phi_i(x)$ of
theorem  1 given by $\phi_i(x)\rightarrow \tilde\phi_i(x)$ where
\begin{eqnarray*}
\Biggl( \tilde\phi_1(x), {\tilde\phi_2(x)\choose\tilde\phi_3(x)},&&
                         {\tilde\phi_4(x)\choose\tilde\phi_5(x)}
\Biggr) =\\ \\
&&  \Biggl( {{ e\sp{\zeta(\nu_2) x}\over-\nu_2} \phi_1(x),
          { \Phi(x;\nu_1)\choose \Phi\sp\prime(x;\nu_1)},
 {e\sp{-\zeta(\nu_2)x} \Phi(x;\nu_2)\choose
   -\nu_2\,e\sp{-\zeta(\nu_2)x} \Phi\sp\prime(x;\nu_2)} \Biggr)}.
\end{eqnarray*}
Now 
$$
\lim_{\nu_2\rightarrow0}\tilde\phi_1(x)=\Phi(x;\nu_1),
$$
and
$$
\lim_{\nu_2\rightarrow0}
{ \biggl|\matrix{\tilde\phi_4(x)&\tilde\phi_4(y)\cr
                 \tilde\phi_5(x)&\tilde\phi_5(y)\cr}\biggr|} 
={\wp }(x)-{\wp }(y) .
$$
Thus (\ref{calfun}) arises as the $\nu_2\rightarrow0$ of (\ref{functional}).

\section{Proof of Theorem 1}
Our proof of Theorem 1 proceeds in two stages. First we will use the
symmetry (\ref{symm}) to transform (\ref{functional}) into a particularly
simple canonical form. This form may be immediately integrated to 
yield a functional
equation studied by Bruschi and Calogero \cite{BCb};  by appealing to
their result our Theorem 1 will follow. The first stage of this process
is entirely algorithmic and consequently we may readily identify the
parameters that appear in our solution. We begin with
\begin{defn}A point $x_0\in\mathbb{C}$ is said to be {\it generic} for
(\ref{functional}) if 
\begin{itemize}
\item[1)] $\phi_k(x)$ is regular at $x_0$ for $k=2\ldots5$,
\item[2)] $\phi_1(x)$ is regular at $2x_0$,
\item[3)]$ { { \biggl|
             \matrix{\phi_2(x_0)&\phi_2\sp\prime(x_0)\cr
                     \phi_3(x_0)&\phi_3\sp\prime(x_0)\cr}\biggr|}
\ne0\quad\quad { \biggl| 
              \matrix{\phi_4(x_0)&\phi_4\sp\prime(x_0)\cr
                      \phi_5(x_0)&\phi_5\sp\prime(x_0)\cr}\biggr|} }\ne0.$
\end{itemize}
\end{defn}

Now let $x_0$ be a generic point. Using at first the matrices $U$ and $V$
of transformation (\ref{symm})
we may choose linear combinations of $\phi_k$ ($k:2\ldots5$) such that
(\ref{functional}) becomes
$$
\tilde\phi_1(x+y)=
{ { \biggl|\matrix{\tilde\phi_2(x)&\tilde\phi_2(y)\cr\tilde\phi_3(x)&
    \tilde\phi_3(y)\cr}\biggr|}
\over
{ \biggl|\matrix{\tilde\phi_4(x)&\tilde\phi_4(y)\cr\tilde\phi_5(x)&
     \tilde\phi_5(y)\cr}\biggr|} },
$$
and (for $k=1,2$)
\begin{equation}
\tilde\phi_{2k}(0)=\tilde\phi_{2k+1}\sp\prime(0)=0,\quad\quad
\tilde\phi_{2k}\sp\prime(0)=\tilde\phi_{2k+1}(0)=1.
\label{init}
\end{equation}
The arguments of the functions have been shifted to be centred on
$x_0$ (or $2x_0$ in the case of $\tilde\phi_1(x)$).
Here we have set
$$
\tilde\phi_1(x)=c\, \phi_1(x+2 x_0)
$$
and (for $k=1,2$)
$$
\Biggl(
\matrix{\tilde\phi_{2k}(x)\cr \tilde\phi_{2k+1}(x)\cr}\Biggr)=
\Biggl(
\matrix{\phi_{2k}\sp\prime(x_0)&\phi_{2k}(x_0)\cr
        \phi_{2k+1}\sp\prime(x_0)&\phi_{2k+1}(x_0)\cr}
\Biggr)\sp{-1}
\Biggl(
\matrix{\phi_{2k}(x+x_0)\cr \phi_{2k+1}(x+x_0)\cr}
\Biggr).
$$
The constant $c$ here is just the ratio of the appropriate determinants
specified in (\ref{symm}).  We next observe

\begin{lem}
For $k=1,2$ we may write
$$
\Biggl(\matrix{\tilde\phi_{2k}(x)\cr \tilde\phi_{2k+1}(x)\cr}\Biggr)=
{1\over\gamma_k(x)}
\Biggl(\matrix{\xi_{k}(x)\cr \xi_{k}\sp\prime(x)\cr}\Biggr),
$$
where $\gamma_k(x), \xi_{k}(x)$ are regular at $0$ and
$$
\xi_{k}(0)=0,\quad\quad \xi_{k}\sp\prime(0)=\gamma_k(0)=1.
$$
Further, upon writing
 $\xi_{k}(x)=e\sp{\lambda_k x}\tilde\xi_{k}(x)$
with $\lambda_k=-\tilde\phi_{2k}\sp{\prime\prime}(0)/2$
the function $\tilde\xi_{k}(x)$, regular at $0$, satisfies
$$
\tilde\xi_{k}(0)=\tilde\xi_{k}\sp{\prime\prime}(0)=0,\quad\quad 
\tilde\xi_{k}\sp\prime(0)=1.
$$

\end{lem}
\begin{proof}
Upon differentiating $\xi_{k}(x)=\gamma_k(x)\tilde\phi_{2k}(x)$ and
comparing with
$\xi_{k}\sp\prime(x)=\gamma_k(x)\tilde\phi_{2k+1}(x)$ we see that
\begin{equation}
{\gamma_k\sp\prime(x)\over \gamma_k(x)}=
{ \tilde\phi_{2k+1}(x)-\tilde\phi_{2k}\sp\prime(x)\over \tilde\phi_{2k}(x)}.
\label{gammak}
\end{equation}
The only issue is whether the righthand side of this differential equation
is regular at $x=0$. Using (\ref{init}) and l'H\^opital's rule
we find 
$$
{\gamma_k\sp\prime(0)\over \gamma_k(0)}=
{ \tilde\phi_{2k+1}\sp\prime(0)-\tilde\phi_{2k}\sp{\prime\prime}(0)
  \over \tilde\phi_{2k}\sp\prime(0)}=-\tilde\phi_{2k}\sp{\prime\prime}(0).
$$
and so $\gamma_k(x)$ and hence $\xi_{k}(x)$ are regular at $0$.
Now $\gamma_k(x)$ is determined by (\ref{gammak}) given an initial
condition which we choose to be $\gamma_k(0)=1$. The remaining
initial conditions for $\xi_{k}(x)$ follow from (\ref{init}).
Indeed from 
$$
\tilde\phi_{2k+1}\sp{\prime}(x)=
{\xi_{k}\sp{\prime\prime}(x)\gamma_k(x)-\xi_{k}\sp{\prime}(x)
\gamma_k\sp\prime(x)\over \gamma_k\sp2(x)}
$$
we also find
$$
\xi_{k}\sp{\prime\prime}(0)=-\tilde\phi_{2k}\sp{\prime\prime}(0).
$$
Now upon writing $\xi_{k}(x)=e\sp{\lambda_k x}\tilde\xi_{k}(x)$
with $\lambda_k=-\tilde\phi_{2k}\sp{\prime\prime}(0)/2$
we obtain the final statement of the lemma.
\end{proof}

Thus far we havent used the exponential part of the symmetry
(\ref{symm}). Utilising this symmetry we set
$\tilde\xi_{0}(x)=e\sp{(\lambda_1-\lambda_2)x}\tilde\phi_1(x)$, 
$\gamma(x)= e\sp{2(\lambda_1-\lambda_2)x}\gamma_2(x)/\gamma_1(x)$.
This scaling has (upon noting $2\lambda_k=\gamma_k\sp\prime(0)$)
the effect of making $\gamma\sp\prime(0)=0$. Thus we obtain
\begin{col}
At any generic point we may rewrite (\ref{functional}) using
the symmetry (\ref{symm}) as 
\begin{equation}
\tilde\xi_{0}(x+y)=\gamma(x)\gamma(y)
{ { \biggl|\matrix{\tilde\xi_1(x)&\tilde\xi_1(y)\cr\tilde\xi_1\sp\prime(x)&
    \tilde\xi_1\sp\prime(y)\cr}\biggr|}
\over
{ \biggl|\matrix{\tilde\xi_2(x)&\tilde\xi_2(y)\cr\tilde\xi_2\sp\prime(x)&
     \tilde\xi_2\sp\prime(y)\cr}\biggr|} },
\label{revfunctional}
\end{equation}
where for $k=1,2$
\begin{equation}
\tilde\xi_{k}(0)=\tilde\xi_{k}\sp{\prime\prime}(0)=\gamma\sp\prime(0)=0,
\quad\quad \tilde\xi_{k}\sp\prime(0)= \gamma(0)=1.
\label{initials}
\end{equation}
\end{col}

Given the complexity of the differential equation (\ref{gammak}) one may
wonder whether (\ref{revfunctional}) simplifies much further. In fact we
find
\begin{lem} 
The functional equation (\ref{revfunctional}), and
consequently (\ref{functional}),  may be written as
\begin{equation}
\partial\biggl({\tilde\xi_1(x+y)\over\tilde\xi_1(x)\tilde\xi_1(y)}\biggl)
=
\partial\biggl({\tilde\xi_2(x+y)\over\tilde\xi_2(x)\tilde\xi_2(y)}\biggl),
\label{newfunctional}
\end{equation}
where $\partial=\partial_x-\partial_y$. Further
\begin{equation}
\tilde\xi_{0}(x)={\tilde\xi_2(x)\over\tilde\xi_1(x)}.
\label{eqnsub}
\end{equation}
\end{lem}
\begin{proof}
Taking the logarithmic derivative of (\ref{revfunctional}) with respect to
$\partial=\partial_x-\partial_y$ we obtain
\begin{equation}
0={\gamma\sp\prime(x)\over\gamma(x)}-{\gamma\sp\prime(y)\over\gamma(y)}
+{\partial\sp2\left(\tilde\xi_1(x)\tilde\xi_1(y)\right)\over\partial
  \left(\tilde\xi_1(x)\tilde\xi_1(y)\right)}
-{\partial\sp2\left(\tilde\xi_2(x)\tilde\xi_2(y)\right)\over\partial
  \left(\tilde\xi_2(x)\tilde\xi_2(y)\right)}.
\label{intermed1}
\end{equation}
Now employing (\ref{initials}) one finds
$$
\partial\left(\tilde\xi_k(x)\tilde\xi_k(y)\right)\vert_{y=0}
=\tilde\xi_k\sp\prime(x)\tilde\xi_k(y)-\tilde\xi_k(x)\tilde\xi_k\sp\prime(y)
\vert_{y=0}=-\tilde\xi_k(x)
$$
and similarly
$$
\partial\sp2\left(\tilde\xi_k(x)\tilde\xi_k(y)\right)\vert_{y=0}=
-2\tilde\xi_k\sp\prime(x).
$$
Upon setting $y=0$ in (\ref{intermed1}) and with these simplifications 
we obtain the differential equation
$$
0={\gamma\sp\prime(x)\over\gamma(x)}+{
2\tilde\xi_1\sp\prime(x)\over\tilde\xi_1(x)}-
{2\tilde\xi_2\sp\prime(x)\over\tilde\xi_2(x)}
$$
with solution
\begin{equation}
\gamma(x)=c\, {\tilde\xi_2\sp2(x)\over\tilde\xi_1\sp2(x)}.
\label{gammadef}
\end{equation}
Again using l'H\^opital's rule and (\ref{initials}) we find the
constant $c=1$.
Therefore (\ref{revfunctional}) may be rewritten as
\begin{equation}
\tilde\xi_{0}(x+y)=
{\tilde\xi_2\sp2(x)\over\tilde\xi_1\sp2(x)}
{\tilde\xi_2\sp2(y)\over\tilde\xi_1\sp2(y)}
{ { \Biggl|\matrix{\tilde\xi_1(x)&\tilde\xi_1(y)\cr\tilde\xi_1\sp\prime(x)&
    \tilde\xi_1\sp\prime(y)\cr}\Biggr|}
\over
{ \Biggl|\matrix{\tilde\xi_2(x)&\tilde\xi_2(y)\cr\tilde\xi_2\sp\prime(x)&
     \tilde\xi_2\sp\prime(y)\cr}\Biggr|} }
=
{ { \Biggl|\matrix{{1\over\tilde\xi_1(x)}&{1\over\tilde\xi_1(y)}\cr
 \big({1\over\tilde\xi_1(x)}\big)\sp\prime&
 \big({1\over\tilde\xi_1(y)}\big)\sp\prime\cr}\Biggr|}
\over
{ \Biggl|\matrix{{1\over\tilde\xi_2(x)}&{1\over\tilde\xi_2(y)}\cr
 \big({1\over\tilde\xi_2(x)}\big)\sp\prime&
 \big({1\over\tilde\xi_2(y)}\big)\sp\prime\cr}\Biggr|}
}
=
{\partial\big({1\over\tilde\xi_1(x)\tilde\xi_1(y)}\big) \over
 \partial\big({1\over\tilde\xi_2(x)\tilde\xi_2(y)}\big)}.
\label{intermed2}
\end{equation}
Letting $y\rightarrow0$ we find
$$
\tilde\xi_{0}(x)={\tilde\xi_2(x)\over\tilde\xi_1(x)}
$$
as required. Utilizing (\ref{eqnsub})  we may immediately rewrite 
(\ref{intermed2})  in the stated form (\ref{newfunctional}).

\end{proof}
We observe that at this stage the symmetry (\ref{symm}) has enabled us
to transform (\ref{functional}) into the form specified by theorem 1.
The solution will follow once we show $1/\tilde\xi_k(x)=\Phi(x;\nu_k)$.
Now (\ref{newfunctional}) may be immediately integrated to give
$$
{\tilde\xi_1(x+y)\over\tilde\xi_1(x)\tilde\xi_1(y)}
={\tilde\xi_2(x+y)\over\tilde\xi_2(x)\tilde\xi_2(y)} +\Theta(x+y).
$$
Upon setting $\alpha(x)={\tilde\xi_2(x)/\tilde\xi_1(x)}$ and
$\psi(x)=\Theta(x)/ \tilde\xi_2(x)$ this may be rearranged into the form
\begin{equation}
{\alpha(x+y)\over\alpha(x)\alpha(y)}=
1+\tilde\xi_2(x)\tilde\xi_2(y)\psi(x+y)
\label{bcalfun}
\end{equation}
which is  the functional equation studied by Bruschi and Calogero
\cite{BCb}.
Calling upon the general analytic solution obtained by these authors,
together with our initial conditions (\ref{initials}), we find\footnote{
For example, from \cite{BCb} we obtain 
$\tilde\xi_2(x)=A e\sp{cx}\sigma(ax\vert\omega,\omega\sp\prime)/
\sigma(ax+\nu\vert\omega,\omega\sp\prime)$. Using the property
$\sigma(ax\vert a\omega,a\omega\sp\prime)=
a\sigma(ax\vert\omega,\omega\sp\prime)$ and the definition of $\Phi(x;\nu)$
we may rewrite this as
$\tilde\xi_2(x)=(A/\sigma(\nu/a)) e\sp{(c-\zeta(\nu/a))x} /\Phi(z;-\nu/a)$.
Now the $x\rightarrow 0$ limit shows 
$(A/\sigma(\nu/a)) e\sp{(c-\zeta(\nu/a))x}=1$.}
$1/\tilde\xi_k(x)=\Phi(x;\nu_k)$ as required.

\section{Proof of Theorem 2}
It is useful at the outset to gather together the various transformations
introduced in the last section:
\begin{eqnarray}
\Biggl(
\matrix{\tilde\phi_{2k}(x)\cr \tilde\phi_{2k+1}(x)\cr}\Biggr)&=&
\Biggl(
\matrix{\phi_{2k}\sp\prime(x_0)&\phi_{2k}(x_0)\cr
        \phi_{2k+1}\sp\prime(x_0)&\phi_{2k+1}(x_0)\cr}
\Biggr)\sp{-1}
\Biggl(
\matrix{\phi_{2k}(x+x_0)\cr \phi_{2k+1}(x+x_0)\cr}
\Biggr) \label{transform1}\\
&=&{1\over\gamma_k(x)}
\Biggl(\matrix{\xi_{k}(x)\cr \xi_{k}\sp\prime(x)\cr}\Biggr) 
={e\sp{\lambda_k x}\over\gamma_k(x)}
\Biggl( \matrix{1&0\cr \lambda_k&1\cr}\Biggr)
\Biggl( \matrix{\tilde\xi_k(x)\cr \tilde\xi_k\sp\prime(x)\cr}\Biggr) 
\label{transform2}\\
&=&{e\sp{\lambda_k x}\over\gamma_k(x)\Phi\sp2(x;\nu_k)}
\Biggl( \matrix{1&0\cr \lambda_k&-1\cr}\Biggr)
\Biggl( \matrix{\Phi(x;\nu_k)\cr \Phi\sp\prime(x;\nu_k)\cr}\Biggr)
\label{transform3}\\
\tilde\xi_0(x)&=&e\sp{(\lambda_1-\lambda_2) x}\,
{\phi_1(x+2 x_0)\over\phi_1(2 x_0)}={\Phi(x;\nu_1)\over\Phi(x;\nu_2)}.
\label{transform4}
\end{eqnarray}
Let us introduce the function
\begin{equation}
f(x)= \gamma_k(x)\Phi\sp2(x;\nu_k)e\sp{-2\lambda_k x}.
\label{transform5}
\end{equation}
Observe that (\ref{gammadef}) entails the function
$f(x)$ is independent of $k$,
$$
\gamma_1(x)\Phi\sp2(x;\nu_1)e\sp{-2\lambda_1 x}
=
\gamma_2(x)\Phi\sp2(x;\nu_k)e\sp{-2\lambda_k x}.
$$
With this definition we may rewrite (\ref{transform1}) and
(\ref{transform3}) to give
\begin{eqnarray}
\left( \matrix{\phi_{2k}(x+x_0)\cr \phi_{2k+1}(x+x_0)\cr} \right)
&=&
{e\sp{-\lambda_k x}\over f(x)}
\Biggl(
\matrix{\phi_{2k}\sp\prime(x_0)&\phi_{2k}(x_0)\cr
        \phi_{2k+1}\sp\prime(x_0)&\phi_{2k+1}(x_0)\cr}
\Biggr)
\Biggl( \matrix{1&0\cr \lambda_k&-1\cr}\Biggr)
\Biggl( \matrix{\Phi(x;\nu_k)\cr \Phi\sp\prime(x;\nu_k)\cr}\Biggr)
\nonumber
\\
\label{transform6}
\\
&=&
{e\sp{-\lambda_k x}\over f(x)}
\left(
\matrix{ \left|
          \matrix{\Phi(x;\nu_k)&\phi_{2k}(x_0)\cr
 \Phi\sp\prime(x;\nu_k)&\phi_{2k}\sp\prime(x_0)+\lambda_k\phi_{2k}(x_0)\cr}
         \right|
        \cr
           \left|
          \matrix{\Phi(x;\nu_k)&\phi_{2k+1}(x_0)\cr
 \Phi\sp\prime(x;\nu_k)&\phi_{2k+1}\sp\prime(x_0)+\lambda_k\phi_{2k+1}(x_0)\cr}
         \right|
        \cr
        }
\right)
\nonumber
\end{eqnarray}
and
\begin{equation}
\Biggl( \matrix{\Phi(x;\nu_k)\cr \Phi\sp\prime(x;\nu_k)\cr}\Biggr)
=
f(x)e\sp{\lambda_k x}\Biggl( \matrix{1&0\cr \lambda_k&-1\cr}\Biggr)
\Biggl(
\matrix{\phi_{2k}\sp\prime(x_0)&\phi_{2k}(x_0)\cr
        \phi_{2k+1}\sp\prime(x_0)&\phi_{2k+1}(x_0)\cr}
\Biggr)\sp{-1}
\Biggl(
\matrix{\phi_{2k}(x+x_0)\cr \phi_{2k+1}(x+x_0)\cr}
\Biggr)
\label{transform7}
\end{equation}
$$
\quad=
{f(x)e\sp{\lambda_k x}\over
\left|
\matrix{\phi_{2k}\sp\prime(x_0)&\phi_{2k}(x_0)\cr
        \phi_{2k+1}\sp\prime(x_0)&\phi_{2k+1}(x_0)\cr}
\right|
 }\,
\left(
\matrix{ \left|
          \matrix{\phi_{2k}(x+x_0)&\phi_{2k}(x_0)\cr
          \phi_{2k+1}(x+x_0)&\phi_{2k+1}(x_0)\cr}
         \right|
        \cr
           \left|
          \matrix{\phi_{2k}(x+x_0)&
                  \phi_{2k}\sp\prime(x_0)+\lambda_k\phi_{2k}(x_0)\cr
          \phi_{2k+1}(x+x_0)&
                  \phi_{2k+1}\sp\prime(x_0)+\lambda_k\phi_{2k+1}(x_0)\cr}
         \right|
        \cr
        }
\right).
\nonumber
$$
Now (\ref{transform4}) and (\ref{transform6}) are of the form stated
in theorem 2 provided we can show $f(x)$, defined in (\ref{transform5}),
can also be put into  the form of the theorem. To see this
note that (\ref{transform1}) shows
\begin{eqnarray}
\tilde\phi_{2k}(x)&=&
{ {\phi_{2k+1}(x_0)\phi_{2k}(x+x_0)- \phi_{2k}(x_0)\phi_{2k+1}(x+x_0)}\over
  {\phi_{2k+1}(x_0)\phi_{2k}\sp\prime(x_0)-
   \phi_{2k}(x_0)\phi_{2k+1}\sp\prime(x_0)}
}
\nonumber
\\
\label{transform8}
\\
&=&
{
\left|
\matrix{\phi_{2k}(x+x_0)&\phi_{2k}(x_0)\cr
        \phi_{2k+1}(x+x_0)&\phi_{2k+1}(x_0)\cr}
\right|
\over
\left|
\matrix{\phi_{2k}\sp\prime(x_0)&\phi_{2k}(x_0)\cr
        \phi_{2k+1}\sp\prime(x_0)&\phi_{2k+1}(x_0)\cr}
\right|
}
\nonumber
\end{eqnarray}
while from (\ref{transform3}) we see
\begin{equation}
\gamma_k(x)={ e\sp{\lambda_k x}\over \Phi(x;\nu_k)\tilde\phi_{2k}(x)}.
\label{transform9}
\end{equation}
%$$
%\gamma_k(x)={e\sp{\lambda_k x}\over \Phi(x;\nu_k)}
%{ {\phi_{2k+1}(x_0)\phi_{2k}\sp\prime(x_0)-
%   \phi_{2k}(x_0)\phi_{2k+1}\sp\prime(x_0)} \over
%  {\phi_{2k+1}(x_0)\phi_{2k}(x+x_0)- \phi_{2k}(x_0)\phi_{2k+1}(x+x_0)}
%}
%$$
Combining these thus shows
\begin{equation}
f(x)=
{e\sp{-\lambda_k x}\over \Phi(x;\nu_k)}\,
{
\left|
\matrix{\phi_{2k}\sp\prime(x_0)&\phi_{2k}(x_0)\cr
        \phi_{2k+1}\sp\prime(x_0)&\phi_{2k+1}(x_0)\cr}
\right|
\over
\left|
\matrix{\phi_{2k}(x+x_0)&\phi_{2k}(x_0)\cr
        \phi_{2k+1}(x+x_0)&\phi_{2k+1}(x_0)\cr}
\right|
}
\label{transform10}
\end{equation}
as required. Also from (\ref{transform8}) and the definition
$\lambda_k=-\tilde\phi_{2k}\sp{\prime\prime}(0)/2$ we find
\begin{eqnarray}
-2\lambda_k &=&
{ {\phi_{2k+1}(x_0)\phi_{2k}\sp{\prime\prime}(x_0)-
   \phi_{2k}(x_0)\phi_{2k+1}\sp{\prime\prime}(x_0)}\over
  {\phi_{2k+1}(x_0)\phi_{2k}\sp\prime(x_0)-
   \phi_{2k}(x_0)\phi_{2k+1}\sp\prime(x_0)}
}
\label{transform11}\\
&=&
\partial_x\ln
{ \Biggl|\matrix{ \phi_{2k}\sp\prime(x+x_0)&\phi_{2k}(x_0)\cr
                       \phi_{2k+1}\sp\prime(x+x_0)&\phi_{2k+1}(x_0)\cr}
  \Biggl| }_{x=0}.
\label{transform12}
\end{eqnarray}

At this stage then we see that if we can determine $\Phi(x;\nu_k)$
all of the terms in (\ref{transform1}-\ref{transform4}) are determined
and we obtain the stated expressions for $\phi_1(x)$, $\phi_2(x)$,
$\phi_3(x)$, $\phi_4(x)$, $\phi_5(x)$ and $f(x)$ given in theorem 2.
It remains therefore to determine the parameters $g_2$, $g_3$
specifying the elliptic functions $\Phi(x;\nu_k)$ as well as
$\nu_1$, $\nu_2$.
To this end we utilise (\ref{transform7}) to give
\begin{eqnarray*}
{\Phi\sp\prime(x;\nu_k) \over \Phi(x;\nu_k)}-\lambda_k
&=&
{
         \left|
          \matrix{\phi_{2k}(x+x_0)& \phi_{2k}\sp\prime(x_0)\cr
          \phi_{2k+1}(x+x_0)& \phi_{2k+1}\sp\prime(x_0)\cr}
         \right|
\over
         \left|
          \matrix{\phi_{2k}(x+x_0)&\phi_{2k}(x_0)\cr
          \phi_{2k+1}(x+x_0)&\phi_{2k+1}(x_0)\cr}
         \right|
}
\\
&=&
\partial_y \ln
         \left|
          \matrix{\phi_{2k}(x+x_0)&\phi_{2k}(y+x_0)\cr
          \phi_{2k+1}(x+x_0)&\phi_{2k+1}(y+x_0)\cr}
         \right|
\Biggl|_{y=0}.
\end{eqnarray*}
Upon using (\ref{defzeta}) to simplify the  left-hand side  of this
equality we thus obtain the first equality of theorem 2,
\begin{equation}
\partial_y \ln
         \left|
          \matrix{\phi_{2k}(x+x_0)&\phi_{2k}(y+x_0)\cr
          \phi_{2k+1}(x+x_0)&\phi_{2k+1}(y+x_0)\cr}
         \right|
\Biggl|_{y=0}
=
\zeta(\nu_k)-\zeta(x)-\zeta(\nu_k-x)-\lambda_k,
\label{transform13}
\end{equation}
and consequently
\begin{equation}
\partial_x\partial_y\ln
{ \Biggl|\matrix{ \phi_{2k}(x+x_0)&\phi_{2k}(y+x_0)\cr
                  \phi_{2k+1}(x+x_0)&\phi_{2k+1}(y+x_0)\cr}
  \Biggr|_{y=0} }
\label{eq:Fdef}
=\wp(x)-\wp(\nu_k-x).
\end{equation}
In fact we have the more general result
\begin{eqnarray*}
&\partial_x\partial_y&\ln
{ \Biggl|\matrix{ \phi_{2k}(x+x_0)&\phi_{2k}(y+x_0)\cr
                  \phi_{2k+1}(x+x_0)&\phi_{2k+1}(y+x_0)\cr}
  \Biggr| }
=
\partial_x\partial_y\ln
{ \Biggl|\matrix{ \Phi(x;\nu_k)&\Phi(y;\nu_k)\cr
                  \Phi\sp\prime(x;\nu_k)&\Phi\sp\prime(y;\nu_k)\cr}
\Biggr| }
\\
&=&
\partial_x\partial_y\ln
\biggl[\Phi(x+y;\nu_k)\Bigl( {\wp }(x)-{\wp }(y) \Bigr)\biggr]
\\
&=&
\wp(x+y)-\wp(\nu_k-x-y)+
{ \wp\sp\prime(x)\wp\sp\prime(y)\over \Bigl( {\wp }(x)-{\wp }(y) \Bigr)\sp2}
\end{eqnarray*}
from which (\ref{eq:Fdef}) arises as the $y\rightarrow0$ limit.

It remains to how that the  Laurent series
of (\ref{transform13}) (or equivalently of (\ref{eq:Fdef})) determines the
parameters of $\Phi(x;\nu_k)$. Set
\begin{equation}
\zeta(\nu_k)-\zeta(x)-\zeta(\nu_k-x)-\lambda_k
=
 -{1\over x}-\lambda_k +\sum_{l=0}F_l\, {x\sp{l+1}\over (l+1)!}
\label{transform15}
\end{equation}
or equivalently
\begin{equation}
\wp(x)-\wp(\nu_k-x)=
{1\over x\sp2}+\sum_{l=0}{F_l\over l!}\, x\sp{l}.
\label{transform16}
\end{equation}
While the coefficients  $F_l$ in these expansions depend on
$k=1,2$ we will avoid including this in our notation: certainly the
combinations of these coefficients giving $g_2$ and $g_3$ 
are independent of $k$. 
Now the left-hand side of (\ref{transform15}) has the expansion
$$
-{1\over x}-\lambda_k -\wp(\nu_k)\, x+\wp\sp{\prime}(\nu_k)\,{x\sp2\over2}+
(2c_2 -\wp\sp{\prime\prime}(\nu_k)){x\sp3\over 3!}+\ldots
$$
while that of (\ref{transform16})  begins
$$
{1\over x\sp2}+c_2\,x\sp2+c_3\,x\sp4+\ldots-\Big\lbrace
 \wp(\nu_k)-x\wp\sp{\prime}(\nu_k)+{x\sp2 \over2}\wp\sp{\prime\prime}(\nu_k)+
 \ldots\Big\rbrace 
$$
From either of these we see
$$
F_0=-\wp(\nu_k),\quad F_1=\wp\sp{\prime}(\nu_k),\quad
F_2=2c_2 -\wp\sp{\prime\prime}(\nu_k),
$$
whereupon utilising  (\ref{eq:wpdef}) we obtain
\begin{equation}
c_2={F_2+6 F_0\sp2 \over 12}={g_2\over20}
\quad{\rm and}\quad
g_3= 6 F_0\sp3 -F_1\sp2 +{5\over3}F_0 F_2.
\label{transform18}
\end{equation}
Thus, as stated in Theorem 2, we may obtain the parameters of the elliptic
functions from the Laurent expansion (\ref{transform15}) for either
choice of $k$, the combinations of the coefficients in (\ref{transform18})
being independent of $k$. The constant terms in the two expansions
the determine $\nu_1$, $\nu_2$ via $F_0=-\wp(\nu_k;g_2,g_3)$.

We now have now established all of Theorem 2. It is perhaps useful to
conclude the section with a lemma that implements the theorem.
\begin{lem}
Let
\begin{eqnarray*}
\partial_x\partial_y\ln
{ \Biggl|\matrix{ \phi_{2k}(x+x_0)&\phi_{2k}(y+x_0)\cr
                  \phi_{2k+1}(x+x_0)&\phi_{2k+1}(y+x_0)\cr}
  \Biggr|_{y=0} }
&=& -{1\over x}-\lambda_k +\sum_{l=0}F_l\, {x\sp{l+1}\over (l+1)!}
\\
&=&
{ h\sp{\prime}(0)h\sp{\prime}(x)\over (h(x)-h(0))\sp2 }
\end{eqnarray*}
where $\displaystyle{ h(x)={\phi_{2k}(x+x_0)\over \phi_{2k+1}(x+x_0)} }$.
Set   $\displaystyle{h_k= {h\sp{(k+1)}(0)\over{ (k+1)!\, h\sp{\prime}(0)}} }$.
Then
$$
F_0= - \left|\matrix{h_1&1\cr h_2&h_1\cr}\right|,\quad
F_1= 2 \left|\matrix{h_1&1&0\cr h_2&h_1&1\cr h_3&h_2&h_1\cr}\right|,\quad
F_2=-6 \left|\matrix{h_1&1&0&0\cr h_2&h_1&1&0\cr h_3&h_2&h_1&1\cr
                     h_4&h_3&h_2&h_1\cr}\right|.
$$
\end{lem}

\section{Examples}
We shall now consider the classical addition theorems of
the Jacobi elliptic functions and then the functional equations
(\ref{dexponential}) and (\ref{eq:biggy}) as examples of our theory.
We have collected several standard results pertaining to
elliptic functions that are of use in our computations in Appendix A.

\noindent{\bf Example 1.} 
As a first application of our theory we consider the addition
theorems for the Jacobi elliptic functions $\dn(x)$, $\cn(x)$ and
$\sn(x)$ where $\dn(x)\equiv \dn(x|m)$ and so on. These may be
cast in the form of (\ref{functional}) as
$$
\dn(x+y)={
 \left|\matrix{\cn\sp\prime(x)&\cn\sp\prime(y) \cr
               \cn(x)&\cn(y)\cr}\right|
         \over
 \left|\matrix{\sn\sp\prime(x)&\sn\sp\prime(y) \cr
               \sn(x)&\sn(y)\cr}\right|
          }\quad{\rm (Jacobi)},\quad\quad
\cn(x+y)={1\over k\sp2}{
 \left|\matrix{\dn\sp\prime(x)&\dn\sp\prime(y) \cr
               \dn(x)&\dn(y)\cr}\right|
         \over
 \left|\matrix{\sn\sp\prime(x)&\sn\sp\prime(y) \cr
               \sn(x)&\sn(y)\cr}\right|
          }
$$
and
$$
\sn(x+y)={
 \left|\matrix{1&1\cr \sn\sp2(x)&\sn\sp2(y)\cr}\right|
         \over
 \left|\matrix{\sn\sp\prime(x)&\sn\sp\prime(y) \cr
               \sn(x)&\sn(y)\cr}\right|
          }\quad{\rm (Cayley)}.
$$

Let us now apply our theorem to the first equality.
The first step is to choose an appropriate generic point $x_0$.
This means we wish $x_0$ to be a regular point for $\cn(x)$, $\sn(x)$
and $\dn(2x)$ as well as
\begin{eqnarray*}
0&\ne& \left|\matrix{\cn\sp\prime(x_0)&\cn\sp{\prime\prime}(x_0)\cr
                   \cn(x_0)&\cn\sp\prime(x_0)\cr}\right|=1-m+m\,\cn\sp4(x_0),
\\
0&\ne& \left|\matrix{\sn\sp\prime(x_0)&\sn\sp{\prime\prime}(x_0)\cr
                   \sn(x_0)&\sn\sp\prime(x_0)\cr}\right|=1-m\,\sn\sp4(x_0).
\end{eqnarray*}
Thus we can take $x_0=0$ for this example.

Using (\ref{transform11}) we find
$$
  -2\lambda_1=\partial_x\ln\cn\sp\prime(x)\big\vert_{x=0} =0\quad{\rm and}\quad
  -2\lambda_2=\partial_x\ln\sn\sp\prime(x)\big\vert_{x=0} =0.
$$
Further, with $h(x)=\phi_2(x)/\phi_3(x)=\cn\sp\prime(x)/\cn(x)$ we obtain
$$
F(x)={1-m+m\,\cn\sp4(x)\over \sn\sp2(x)\dn\sp2(x)}
    ={1\over x\sp2}+ {1-2 m\over3}+{1+14 m-14 m\sp2 \over 15}x\sp2+\ldots
$$
while with $h(x)=\phi_4(x)/\phi_5(x)=\sn\sp\prime(x)/\sn(x)$ we obtain
$$
F(x)={1-m\,\sn\sp4(x)\over \sn\sp2(x)}
    ={1\over x\sp2}+ {1+m\over3}+{1-16 m+m\sp2 \over 15}x\sp2+\ldots
$$
In both cases we find
$$
g_2={4\over3}(1-m+m\sp2)\quad{\rm and}\quad
g_3={4\over27}(m-2)(2m-1)(m+1),
$$ 
(the required equality providing a nontrivial check) which means
$$
e_1={2-m\over3},\quad e_2={2 m-1\over3}\quad{\rm and}\quad e_3={-1-m\over3}.
$$
Further,
$$\wp(\nu_1)={2 m-1\over3}\quad{\rm and}\quad
  \wp(\nu_2)={-1-m\over3}.
$$
Comparison with (\ref{ehoms}) and (\ref{jacgs})
shows $\omega=K(m)$, $\omega\sp\prime=iK\sp\prime(m)$,
$\nu_1=K(m)+iK\sp\prime(m)$ and $\nu_2=K\sp\prime(m)$.
We may also calculate $f(x)=1/\sn\sp2(x)$ and upon
using (\ref{eq:phiJacs}) our identity may be rewritten as
\begin{eqnarray*}
\dn(x+y)&=&{\Phi(x+y;K(m)+iK\sp\prime(m))\over \Phi(x+y;iK\sp\prime(m))} \\
&=&
{{\biggl| \matrix{\Phi(x;K(m)+iK\sp\prime(m))&\Phi(y;K(m)+iK\sp\prime(m))\cr
            \Phi(x;K(m)+iK\sp\prime(m))\sp\prime&
            \Phi(y;K(m)+iK\sp\prime(m))\sp\prime\cr}\biggr| }\over
 {\biggl| \matrix{\Phi(x;iK\sp\prime(m))&\Phi(y;iK\sp\prime(m))\cr
            \Phi(x;iK\sp\prime(m))\sp\prime&
            \Phi(y;iK\sp\prime(m))\sp\prime\cr}\biggr| }
}\\
&=&
{{1\over\sn\sp2(x)}
 \left|\matrix{\cn\sp\prime(x)&\cn\sp\prime(y) \cr
               \cn(x)&\cn(y)\cr}\right|
         \over
 {1\over\sn\sp2(x)}
 \left|\matrix{\sn\sp\prime(x)&\sn\sp\prime(y) \cr
               \sn(x)&\sn(y)\cr}\right|
}\\
\end{eqnarray*}

The second identity may be treated in the same manner yielding
$\nu_1=K(m)$ and $\nu_2=K\sp\prime(m)$.
The third identity is a little different. It may be rewritten as
\begin{eqnarray*}
\Phi(x+y;iK\sp\prime(m))&=&
{1\over \sn(x+y)}=
{ {1\over\sn\sp2(x)} 
  \left|\matrix{\sn\sp\prime(x)&\sn\sp\prime(y) \cr
               \sn(x)&\sn(y)\cr}\right|
   \over
  {1\over\sn\sp2(x)}
  \left|\matrix{1&1\cr \sn\sp2(x)&\sn\sp2(y)\cr}\right|
}\\
&=&
{ {\biggl| \matrix{\Phi(x;iK\sp\prime(m))&\Phi(y;iK\sp\prime(m))\cr
            \Phi(x;iK\sp\prime(m))\sp\prime&
            \Phi(y;iK\sp\prime(m))\sp\prime\cr}\biggr| }
  \over 
  {\biggl| \matrix{\Phi\sp2(x;iK\sp\prime(m))&\Phi\sp2(y;iK\sp\prime(m))\cr
            1&1\cr}\biggr| }
}.
\end{eqnarray*}
Now
$$
\Phi\sp2(x;iK\sp\prime(m))-\Phi\sp2(y;iK\sp\prime(m))
=\wp(x)-\wp(y)
$$
and the  required identity follows from the general solution by the limiting
procedure described in section 3.

\noindent{\bf Example 2.}
We shall now determine the general analytic solution of
$$
\phi_1(x+y)= \phi_4(x)\phi_5(y)+\phi_4(y)\phi_5(x)=
{ { \biggl|\matrix{\phi_2(x)&\phi_2(y)\cr\phi_3(x)&\phi_3(y)\cr}\biggr|}
\over
{ \biggl|\matrix{\phi_4(x)&\phi_4(y)\cr\phi_5(x)&\phi_5(y)\cr}\biggr|} },
$$
where $ \phi_2(x)=\phi_4\sp2(x)$, $\phi_3(x)=\phi_5\sp2(x)$.
The particular case $\phi_1(x)=\phi_4(x)$ was treated in
\cite{BK}.

Suppose $x_0$ is a generic point. Then from
$$
0\ne 
\biggl|
       \matrix{\phi_2(x_0)&\phi_2\sp\prime(x_0)\cr
               \phi_3(x_0)&\phi_3\sp\prime(x_0)\cr}
\biggr|
=
2 \phi_4(x_0) \phi_5(x_0)
\biggl|
       \matrix{\phi_4(x_0)&\phi_4\sp\prime(x_0)\cr
               \phi_5(x_0)&\phi_5\sp\prime(x_0)\cr}
\biggr|
$$
we see $\phi_4(x_0)\ne0$, $\phi_5(x_0)\ne0$ and 
$\phi_1(2 x_0)= 2 \phi_4(x_0) \phi_5(x_0)\ne0$.
Further, from (\ref{transform11}), we find
\begin{equation}
\lambda_1=\lambda_2-{1\over2}
\left( {\phi_4\sp\prime(x_0)\over \phi_4(x_0)}
     + {\phi_5\sp\prime(x_0)\over \phi_5(x_0)}
\right).
\label{eq:lams}
\end{equation}

Our strategy  is as follows. We will first determine
$\nu_1$, $\nu_2$, $\lambda_1$, $\lambda_2$, the parameters
describing the elliptic functions and  the ratio 
$\phi_4(x+x_0)\phi_5(x_0) / \phi_5(x+x_0)\phi_4(x_0)$.
Then from
\begin{eqnarray*}
\phi_1(x+2 x_0) 
&=& \phi_4(x+x_0)\phi_5(x_0)+\phi_4(x_0)\phi_5(x+x_0) \\
&=&\phi_4(x+x_0)\phi_5(x_0)\left(
 1+{\phi_5(x+x_0)\phi_4(x_0)\over \phi_4(x+x_0)\phi_5(x_0)}
 \right) \\ 
&=&
e\sp{(\lambda_2-\lambda_1) x}\,{\Phi(x;\nu_1)\over\Phi(x;\nu_2)}
\phi_1(2 x_0)
\end{eqnarray*}
we will obtain
\begin{equation}
\phi_4(x+x_0)=
{2 \phi_4(x_0) e\sp{(\lambda_2-\lambda_1) x} \over
 1+{\phi_5(x+x_0)\phi_4(x_0) / \phi_4(x+x_0)\phi_5(x_0) }
} 
{\Phi(x;\nu_1)\over\Phi(x;\nu_2)}
\label{eq:strategy}
\end{equation}
with $\phi_5(x+x_0)$ immediately following.

Now from (\ref{transform6}) we obtain 
\begin{equation}
{\phi_{2k}(x+x_0)\over \phi_{2k+1}(x+x_0)}
{\phi_{2k+1}(x_0)\over \phi_{2k}(x_0)}
=
1+ {N_k\over D_k}
\label{eq:NDks}
\end{equation}
%$$
%{\phi_{2k}(x+x_0)\over \phi_{2k+1}(x+x_0)}
%{\phi_{2k+1}(x_0)\over \phi_{2k}(x_0)}
%=
%1+
%{ \phi_{2k}\sp\prime(x_0)/\phi_{2k}(x_0)
%  - \phi_{2k+1}\sp\prime(x_0)/\phi_{2k+1}(x_0)
% \over
%  \phi_{2k+1}\sp\prime(x_0)/\phi_{2k+1}(x_0)+\lambda_k
%  - \Phi\sp\prime(x;\nu_k) / \Phi(x;\nu_k)
%}
%$$
where
$$
N_k= { \phi_{2k}\sp\prime(x_0)\over\phi_{2k}(x_0)}
    -{ \phi_{2k+1}\sp\prime(x_0)\over\phi_{2k+1}(x_0)}
\quad\quad
D_k= { \phi_{2k+1}\sp\prime(x_0)\over\phi_{2k+1}(x_0)}+\lambda_k
    -{\Phi\sp\prime(x;\nu_k)\over\Phi(x;\nu_k)}
$$
Here $N_1=2 N_2$ and by our assumption that $x_0$ was a generic point
these are non vanishing. Further, from 
$ \phi_2(x)=\phi_4\sp2(x)$ and $\phi_3(x)=\phi_5\sp2(x)$, we see that
$$
1+N_1/D_1 =\left( 1+ N_2/D_2\right)\sp2.
$$
Expanding this shows $D_2\sp2=(D_2+N_2/2)D_1$ which upon using
(\ref{eq:lams}) yields
\begin{eqnarray}
\left(
 { \phi_5\sp\prime(x_0)\over\phi_5(x_0)} +\lambda_2
 -{\Phi\sp\prime(x;\nu_2)\over\Phi(x;\nu_2)}
\right)\sp2
&=&
\left(
 { \phi_5\sp\prime(x_0)\over\phi_5(x_0)} +\lambda_2
 +{N_2\over2}
%\left( { \phi_4\sp\prime(x_0)\over\phi_4(x_0)}-
%       { \phi_5\sp\prime(x_0)\over\phi_5(x_0)}
% \right)
 -{\Phi\sp\prime(x;\nu_2)\over\Phi(x;\nu_2)}
\right)
\\ \nonumber &&\quad
\left(
 { \phi_5\sp\prime(x_0)\over\phi_5(x_0)} +\lambda_2
 -{N_2\over2}
%\left( { \phi_4\sp\prime(x_0)\over\phi_4(x_0)}-
%       { \phi_5\sp\prime(x_0)\over\phi_5(x_0)}
% \right)
 -{\Phi\sp\prime(x;\nu_1)\over\Phi(x;\nu_1)}
\right)
\label{eq:poles}
\end{eqnarray}
Suppose that $\nu_2$ is finite. Comparing the pole behaviour of each
side of (\ref{eq:poles}) shows that $\nu_1=\nu_2$ and  consequently
that $N_2=0$, a contradiction.
The remaining  possibility is that $\nu_2$ is infinite which
we now show to be a consistent solution. This can only
happen if the elliptic function degenerates into a hyperbolic or
trigonometric function and without loss of generality we choose the former.
In this case
\begin{equation}
\Phi(x;\nu)={\kappa\sinh \kappa(\nu-x)\over\sinh \kappa\nu
               \sinh \kappa x}
e\sp{ x \kappa\coth\kappa\nu}\quad{\rm and}\quad
\Phi(x;\infty)={\kappa\over\sinh \kappa x}.
\label{eq:phiinff}
\end{equation}
Let us the suppose $\nu_2=\infty$. 
Utilising (\ref{Phiinf}) and (\ref{Phiinfs}) we then  must solve
$$
\left(
 { \phi_5\sp\prime(x_0)\over\phi_5(x_0)} +\lambda_2
 +\kappa\, \coth\kappa x
\right)\sp2
=
\left(
 { \phi_5\sp\prime(x_0)\over\phi_5(x_0)} +\lambda_2
 +{N_2\over2}
 +\kappa\, \coth\kappa x
\right)\quad\quad\quad
$$
$$
\left(
 { \phi_5\sp\prime(x_0)\over\phi_5(x_0)} +\lambda_2
 -{N_2\over2}
 +\kappa\coth\kappa (\nu_1-x) +\kappa\coth\kappa x-\kappa\coth\kappa\nu_1
\right).
$$
This holds provided
$$
N_2\sp2= {4 \kappa\sp2\over\sinh \kappa\nu_1},
\quad\quad
{ \phi_5\sp\prime(x_0)\over\phi_5(x_0)} +\lambda_2
 +{N_2\over2}
 +\kappa\, \coth\kappa \nu_1=0
$$
which determines $\nu_1$, $\lambda_2$ and (via (\ref{eq:lams}))
$\lambda_1$ in terms of $\phi_{4}(x_0)$, $\phi_{5}(x_0)$,
$\phi_4\sp\prime(x_0)$  and $\phi_5\sp\prime(x_0)$.
The choice of sign in taking the square-root here is
arbitrary (just defining $\nu_1$) and we will take
$N_2 = -2\kappa/\sinh \kappa\nu_1$.
Substituting these into (\ref{eq:NDks}) we find
\begin{eqnarray*}
{\phi_{4}(x+x_0)\over \phi_{5}(x+x_0)}
{\phi_{5}(x_0)\over \phi_{4}(x_0)}
&=&
{\kappa\coth\kappa x- \kappa\, \coth\kappa \nu_1+{N_2/2}
\over \kappa\coth\kappa x- \kappa\, \coth\kappa \nu_1-{N_2/2}
}
\\ &=&
\coth(\kappa \nu_1/2) \tanh\kappa( \nu_1/2-x).
\end{eqnarray*}
Now employing (\ref{eq:phiinff}) shows
\begin{equation}
\phi_1(x+2 x_0)=e\sp{(\lambda_2-\lambda_1+\kappa\coth\kappa\nu_1) x}\,
{\sinh\kappa(\nu_1-x)\over\sinh\kappa \nu_1} \, \phi_1(2 x_0)
\label{sol1}
\end{equation}
where the exponential may be rewritten to yield
\begin{eqnarray*}
\lambda_2-\lambda_1+\kappa\coth\kappa\nu_1 &=&
{\phi_4\sp\prime(x_0)\over \phi_4(x_0)} -{N_2\over2}+\kappa\coth\kappa\nu_1=
{\phi_4\sp\prime(x_0)\over \phi_4(x_0)}+\kappa\coth(\kappa\nu_1/2)\\
&=&{\phi_5\sp\prime(x_0)\over \phi_5(x_0)} +{N_2\over2}+\kappa\coth\kappa\nu_1
={\phi_5\sp\prime(x_0)\over \phi_5(x_0)} +\kappa\tanh(\kappa\nu_1/2).
\end{eqnarray*}

We now have the information needed to determine $\phi_{4}(x)$ and
$\phi_{5}(x)$ via (\ref{eq:strategy}) which gives
\begin{eqnarray}
\phi_{4}(x+x_0)&=&{\sinh\kappa( \nu_1/2-x)\over \sinh(\kappa \nu_1/2)}
e\sp{\left(\phi_4\sp\prime(x_0)/ \phi_4(x_0)+\kappa\coth(\kappa\nu_1/2)
     \right)x}\, \phi_{4}(x_0)
\label{sol2}\\
\phi_{5}(x+x_0)&=&{\cosh\kappa( \nu_1/2-x)\over \cosh(\kappa \nu_1/2)}
e\sp{\left(\phi_5\sp\prime(x_0)/ \phi_5(x_0)+\kappa\tanh(\kappa\nu_1/2)
     \right)x}\, \phi_{5}(x_0)
\label{sol3}
\end{eqnarray}
Assembling this provides
\begin{lem}
The general analytic solution to (\ref{dexponential}) is given by
(\ref{sol1}), (\ref{sol2}) and (\ref{sol3}),  where  $x_0$ is a generic 
point.
\end{lem}

\noindent{\bf Example 3.} We conclude by showing how our results determine
the solutions of the functional equation (\ref{eq:biggy}):
\[
\Psi _1(x+y)=\Psi _2(x+y) \phi_2(x)\phi_3(y) +\Psi_3(x+y) \phi_4(x)\phi_5(y).
\]
This equation encompasses as particular cases the equations
(\ref{dexponential}) (with $\Psi _2=\Psi _3=1$, $\phi _2(x)=\phi _4(x)$
and $\phi _3(x)=\phi _5(x)$) and (\ref{calogero}) (with
$(\phi _2(x)=\phi_3(x)$ and $\phi _4(x)=\phi _5(x)$) which have
already been discussed. Because of
this we will only consider the generic case
$\phi_2(x)\neq \lambda \phi_3(x)$, $\phi_4(x)\neq \gamma \phi_5(x)$
and $\Psi_2(x)\neq \delta \Psi_3(x)$
(where $\lambda ,\gamma ,\delta $ are constants) rather than these limits.
Our first step is to relate (\ref{eq:biggy}) to (\ref{functional}):

\begin{lem}
The functions $\Psi_m(x)$ ($m=1,2,3$) and $\phi_n(x)$ ($n=2,3,4,5$) give
a solution of equation (\ref{eq:biggy}) if and only if 
\begin{equation}
\frac{\Psi_3(x+y)}{\Psi_2(x+y)}= -
{ { \biggl|\matrix{\phi_2(x)&\phi_2(y)\cr\phi_3(x)&\phi_3(y)\cr}\biggr|}
\over
{ \biggl|\matrix{\phi_4(x)&\phi_4(y)\cr\phi_5(x)&\phi_5(y)\cr}\biggr|} }
\label{eq:inter1}
\end{equation}
and
\begin{equation}
\frac{\Psi_1(x+y)}{\Psi _2(x+y)}=-
{ { \biggl|\matrix{\phi_2(x) \phi_5(x)&\phi_2(y)\phi_5(y)\cr
                   \phi_3(x) \phi_4(x)&\phi_3(y)\phi_4(y)\cr}\biggr|}
\over
{ \biggl|\matrix{\phi_4(x)&\phi_4(y)\cr\phi_5(x)&\phi_5(y)\cr}\biggr|} }
\label{eq:inter2}
\end{equation}
\end{lem}

\begin{proof} Assume first that the functions ${\Psi_m, \phi_n}$
give a solution of equation (\ref{eq:biggy}). Then after
interchanging $x$ and $y$ in (\ref{eq:biggy}) and subtracting
the result from (\ref{eq:biggy}) we obtain equation (\ref{eq:inter1}). 
Upon substituting the formula for  ${\Psi_3(x+y)}/{\Psi_2(x+y)}$
into (\ref{eq:biggy}) we arrive the formula (\ref{eq:inter2}).

In the other direction, let the  functions ${\Psi_m, \phi_n}$
now satisfy (\ref{eq:inter1}) and (\ref{eq:inter2}). Upon writing the
right hand side of (\ref{eq:biggy}) as
\begin{equation}
\begin{array}{l}
\Psi_2(x+y)\phi_2(x)\phi_3(y)+\Psi_3(x+y)\phi_4(x)\phi_5(y)=\\
\quad\quad\quad\Psi_2(x+y)\left( \phi_2(x)\phi_3(y) 
 +\frac{\Psi_3(x+y)}{\Psi_2(x+y)}\phi_4(x)\phi_5(y)\right)
\end{array}
\label{eq:inter3}
\end{equation}
and using expression  (\ref{eq:inter1}) for 
${\Psi_3(x+y)}/{\Psi_2(x+y)}$ we find the term in brackets  
in (\ref{eq:inter3})  rearranges to give precisely the
right hand side of (\ref{eq:inter2}); substituting for this
then yields (\ref{eq:biggy}) and therefore the required solution.
\end{proof}

We may now apply theorem 1 to show that if $\Psi_m(x)$ 
($m=1,2,3$) give a solution of (\ref{eq:biggy}) then we must have the
ratios
\begin{equation}
\frac{\Psi_1(x)}{\Psi_2(x)}= c_1e^{\lambda_1x}
\frac{\Phi(x;\mu_1)}{\Phi(x;\mu_2)},\quad
\frac{\Psi_3(x)}{\Psi_2(x)}=c_2e^{\lambda_2x}
\frac{\Phi(x;\mu_3)}{\Phi(x;\mu_4)}.
\label{eq:biggyratios}
\end{equation}
Further, because the denominators of (\ref{eq:inter1}) and (\ref{eq:inter2})
are the same, theorem 2 shows that $\mu_4=\mu_2$.  Theorem 1 also
determines the functions $\phi_n(x)$ ($n=2,3,4,5$) up to a ${\cal G}$ action.
In fact, given three functions $\Psi_m(x)$ 
($m=1,2,3$) whose ratios satisfy (\ref{eq:biggyratios}) with $\mu_4=\mu_2$,
this is also sufficient to guarantee there are functions 
$\phi_n(x)$ ($n=2,3,4,5$) for which (\ref{eq:biggy}) holds true.
To see this  let us substitute these ratios into equation (\ref{eq:biggy}) 
to give 
\begin{eqnarray}
c_1e^{\lambda_1 (x+y)}\Phi(x+y;\mu_1)=\Phi(x+y;\mu _2)
\phi_2(x)\phi_3(y)+\nonumber \\
c_2e^{\lambda_2 (x+y)}\Phi(x+y;\mu_3)\phi_4(x)\phi_5(y).
\label{eq:biggytr}
\end{eqnarray}

We will have established sufficiency
once we have shown  how to construct the functions 
$\phi_n(x)$.This will be achieved 
utilizing various properties of the functions $\Phi (x;\nu )$.
\begin{lem}
The Baker-Akhiezer functions $\Phi (x;\nu )$  satisfy the equations
\begin{equation}
\Phi (x+\alpha ;\nu )=-e^{(\zeta(\alpha -\nu )+\zeta (\nu )-\zeta(\alpha ))x}
\Phi (\alpha ;\nu )\
\frac{ \Phi (x;\nu -\alpha )}{\Phi (-x;\alpha )}
\label{eq:biggypf2}
\end{equation}
and
\begin{eqnarray}
c e^{\gamma (x+y)}\,\Phi(x+y;\nu_1+\nu_2)=\Phi(x+y;\nu_1)\,\Phi(x;\nu_2)
\,\Phi (y;\nu_2)-\nonumber \\
\Phi (x+y;\nu_2)\,\Phi(x;\nu_1)\,\Phi(y;\nu_1),
\label{eq:biggylem}
\end{eqnarray}
where $c=\wp(\nu _2)-\wp(\nu _1)$ and 
$\gamma =\zeta (\nu _1)+\zeta (\nu _2)-\zeta (\nu _1+\nu _2)$.
\end{lem}
These follow directly from the definition of $\Phi (x;\nu )$  and properties
of the Weierstrass sigma function; in particular
(\ref{eq:biggylem}) is a consequence of the \lq three term relation\rq\
of the sigma function (\cite[20.53, {\it Ex:5}]{WW}).

%\begin{proof}We will use the same method of proof as applied in (\cite{BK}) 
%to establish the addition theorems for Baker-Akhiezer functions.
%The formula 
%\[
%\tilde\Phi(x):=
%\Phi (x+\alpha ;\nu )=\frac{\sigma (\nu -\alpha- x)}{\sigma (\nu
%)\sigma (x+\alpha )}e\sp{\zeta (\nu )(x+\alpha )} 
%\]
%defines a unique function $\tilde\Phi (x)$ such that
%\begin{enumerate}
%\item  $\tilde\Phi $ is a meromorphic function with simple pole at the point 
%$-\alpha $;
%\item
%$\tilde\Phi $ possesses the periodicity properties:
%\[
%\tilde
%\Phi(x+\omega _\ell )=e^{\zeta (\nu )\omega _\ell -\eta _\ell \nu }
%\tilde\Phi (x), 
%\]
%where $\omega _\ell $ are the periods.
%\item $\tilde\Phi $ is normalized by the condition 
%$\tilde\Phi (0)= \Phi (\alpha ;\nu )$.
%\end{enumerate}

%To prove the statement of the Lemma it is enough to check that the
%right-hand side part of (\ref{eq:biggylem})
%has the same analytical properties as the
%function $\tilde\Phi (x)$ that uniquely defines the function $ce^{\gamma
%(x+y)}\Phi (x+y;\nu _1+\nu _2)$.
%\end{proof}

Upon setting $x\rightarrow x+\alpha$ in (\ref{eq:biggylem})
we obtain 
\begin{equation}
\begin{array}{rl}
ce^{\gamma (x+y+\alpha)}\Phi (x+y+\alpha;\nu _1+\nu _2)&=\Phi(x+y+\alpha ;\nu_1)
\Phi (x+\alpha ;\nu_2)\Phi (y;\nu_2)\nonumber \\
&\quad-\Phi (x+y+\alpha ;\nu _2)\Phi (x+\alpha ;\nu _1)\Phi (y;\nu _1)
\label{eq:biggypf1}
\end{array}
\end{equation}

Now by substituting (\ref{eq:biggypf2}) in (\ref{eq:biggypf1})
and setting $\mu _1=\nu _1+\nu _2-\alpha$, $\mu_2=\nu_1-\alpha$ and
$\mu _3=\nu_2-\alpha$ one obtains after some rearrangement
\begin{equation}
\begin{array}{rl}
c\sp\prime e^{\lambda\sp\prime (x+y)}\Phi (x+y;\mu_1)=&
\Phi (x+y;\mu_2) \frac{\Phi (x;\mu_3)}{\Phi (-x;\mu_1-\mu_2-\mu_3)}
\Phi (y;\mu_1-\mu_2)+\nonumber\\
& c\sp{\prime\prime} e^{\lambda\sp{\prime\prime} (x+y)}
\Phi (x+y;\mu_3) \frac{\Phi (x;\mu_2)}{\Phi (-x;\mu_1-\mu_2-\mu_3)}
\Phi (y;\mu_1-\mu_3)
\label{eq:biggypf3}
\end{array}
\end{equation}
for appropriate constants $c\sp\prime, c\sp{\prime\prime},
\lambda\sp\prime, \lambda\sp{\prime\prime} $.
This is precisely of the desired form (\ref{eq:biggytr}).
Therefore we have shown
\begin{thm}
Given functions $\Psi_m(x)$ ($m=1,2,3$), there are functions $\phi_n(x)$
($n=2,3,4,5$) for which the functional equation
(\ref{eq:biggy}) is true if and only if the following ratios take place:
\begin{equation}
\frac{\Psi_1(x)}{\Psi_2(x)}= c_1e^{\lambda_1x}
\frac{\Phi(x;\mu_1)}{\Phi(x;\mu_2)},\quad
\frac{\Psi_3(x)}{\Psi_2(x)}=c_2e^{\lambda_2x}
\frac{\Phi(x;\mu_3)}{\Phi(x;\mu_2)},
\label{eq:biggythm}
\end{equation}
where $c_m,\lambda_m,$ ($m=1,2$) and $\mu_n$, ($n=1,2,3$) are free parameters.
\end{thm}

\vskip 1in

\noindent{\bf Acknowledgements:} 
One of us (V.M.B.) thanks the Royal Society for a Kapitza Fellowship
in $1993$ and the EPSRC for a Visiting Fellowship in $1994-5$ over which time
this paper developed.

\newpage
\appendix
\section{Elliptic Functions}
\subsection{ The Weierstrass elliptic functions:}
The Weierstrass elliptic functions are based on a lattice with periods
$2\omega$ and $2\omega\sp\prime$,  where $\Im(\omega\sp\prime/\omega)>0$.
They satisfy the homogeneity relations
\begin{eqnarray}
\sigma(tx|t\omega,t\omega\sp\prime)&=&t\sigma(x|\omega,\omega\sp\prime),
\quad
\zeta(tx|t\omega,t\omega\sp\prime)=t\sp{-1}\zeta(x|\omega,\omega\sp\prime),
\\
\wp(tx|t\omega,t\omega\sp\prime)&=&t\sp{-2}\wp(x|\omega,\omega\sp\prime).
\end{eqnarray}
Here
$\zeta(x)=\left(\ln\sigma(x)\right)\sp\prime$ and
$\wp(x)=-\zeta\sp\prime(x)$.
These homogeneity relations
mean our function $\Phi(x;\nu)\equiv\Phi(x;\nu|\omega,\omega\sp\prime)$
satisfies
$$
\Phi(t x;t\nu|t\omega,t\omega\sp\prime)=
t\sp{-1}\Phi(x;\nu|\omega,\omega\sp\prime).
$$

The parameters $g_2$ and $g_3$ alternately used to describe the elliptic
function are given by
$$
g_2=60\mathop{{\sum}'}_{m,n\in\mathbb{Z}} \Omega\sp{-4},\quad\quad
g_3=140\mathop{{\sum}'}_{m,n\in\mathbb{Z}} \Omega\sp{-6},
$$
where $\Omega=2 m\omega+2 n \omega\sp\prime$.

The Weierstrass $\wp$-function, $\wp(x)\equiv\wp(x|\omega,\omega\sp\prime)=
\wp(x;g_2,g_3)$, satisfies the differential equation
\begin{equation}
\wp\sp{\prime2}(x)=4 \wp\sp{3}(x)-g_2 \wp(x)- g_3=
4\Big( \wp\sp{3}(x)-5 c_2\wp(x)-7 c_3\Big)
\end{equation}
where $c_2={g_2\over20}$ and $c_3={g_3\over28}$. This means
\begin{equation}
\wp\sp{\prime\prime}(x)=6\wp\sp{2}(x)-10 c_2.
\label{eq:wpdef}
\end{equation}
The   terms $c_k$ ($k\ge4$) in the Laurent expansion of the Weierstrass
$\wp$ function,
\begin{equation}
\wp(x)={1\over x\sp2}+\sum_{l=2}c_l\, x\sp{2l-2},
\label{wplaurent}
\end{equation}
are expressible in terms of $c_2$ and $c_3$.

If $e_i\equiv e_i(\omega,\omega\sp\prime)$ ($i=1,2,3$)
denote the roots of the cubic
\begin{equation}
4 x^3-g_2 x-g_3=0,
\label{eq:cubic}
\end{equation}
whence 
$$\wp\sp{\prime2}(x)=4(\wp(x)-e_1)(\wp(x)-e_2)(\wp(x)-e_3),$$
the half-periods $\omega_i$ are defined by 
$$\wp(\omega_i)=e_i,\quad
{\rm where}\quad
\omega_1=\omega,\quad\omega_2=\omega+\omega\sp\prime,\quad{\rm and}\quad
\omega_3=\omega\sp\prime.
$$
Clearly
\begin{equation}
e_i(t\omega,t\omega\sp\prime)=t\sp{-2}e_i(\omega,\omega\sp\prime).
\label{ehoms}
\end{equation}

Assuming that $g_2$ and $g_3$ are real, and that the discriminant
of (\ref{eq:cubic}) is positive (the case of interest in this paper)
the $e_i$'s are real and may be ordered $e_1\geq e_2\geq e_3$.

\subsection{The Jacobi elliptic functions:}
The Jacobi elliptic functions are characterised by a parameter $m$.
Thus  for example $\sn(x)\equiv\sn(x|m)$ has  periods $4K(m)$
and $2 i K\sp\prime(m)$, where $K(m)$ is the complete elliptic
function of the first kind.
The $e_i$'s are related to the parameter  $m$ of the Jacobi elliptic
functions by
\begin{equation}
e_1={2-m\over3}{K\sp2(m)\over\omega\sp2},
\quad e_2={2 m-1\over3}{K\sp2(m)\over\omega\sp2}
\quad{\rm and}\quad e_3={-1-m\over3}{K\sp2(m)\over\omega\sp2}.
\label{weierjac}
\end{equation}
Thus
\begin{equation}
g_2={4\over3}(1-m+m\sp2){K\sp4(m)\over\omega\sp4},\quad\quad
g_3={4\over27}(m-2)(2m-1)(m+1){K\sp6(m)\over\omega\sp6}.
\label{jacgs}
\end{equation}
and
\begin{equation}
{\omega\over \omega\sp\prime}={iK\sp\prime(m)\over K(m)},\quad\quad
\omega={K(m)\over\sqrt{e_1-e_3}}.
\end{equation}
We find
$$
\Phi(x;\omega\sp\prime)={a\over\sn(a\, x)},\quad
\Phi(x;\omega)=
  {a}\,{\cn(a\, x)\over\sn(a\, x)}
\quad{\rm and}\quad
\Phi(x;\omega+\omega\sp\prime)=
  {a}\,{\dn(a\, x)\over\sn(a\, x)}.
$$
Here $a=\sqrt{e_1-e_3}$ converts the periods of $\Phi$ based on
$\omega$, $\omega\sp\prime$ to those of the Jacobi functions based on
$K(m)$, $iK\sp\prime(m)$. Equally we may write this as
\begin{equation}
\Phi(x;t\omega\sp\prime|t\omega, t\omega\sp\prime)=
{1\over\sn(x|m)}, \quad{\rm with}\quad t=\sqrt{e_1-e_3}={K(m)\over\omega}.
\end{equation}
Thus, with these periods in mind, we write
\begin{eqnarray}
\Phi(x;iK\sp\prime(m))&=&{1\over\sn(x)},\quad
\Phi(x;K(m)+iK\sp\prime(m))={\dn( x)\over\sn(x)}\\
\Phi(x;K(m))&=&{\cn( x)\over\sn(x)}
\label{eq:phiJacs}
\end{eqnarray}

\subsection{Degenerations:}
When the discriminant of (\ref{eq:cubic}) vanishes one (or both) of the
periods of the elliptic function vanishes yielding hyperbolic, trigonometric
(or rational) functions. If $e_1=e_2=c$, $e_3=-2c$ (and so
$g_2=12c\sp2$, $g_3 = -8c\sp3$) we then have
$$
\sigma(x;12c\sp2, -8c\sp3)={\sinh \kappa x\over \kappa}
e\sp{-\kappa\sp2 x^2/6}
\quad{\rm and}\quad
\wp(x;12c\sp2, -8c\sp3)={\kappa\sp2\over3}+{\kappa\sp2\over\sinh\sp2\kappa x}
$$
where $\kappa=\sqrt{3c}$.
In this case
\begin{eqnarray}
\Phi(x;\nu)&=& {\kappa\sinh \kappa(\nu-x)\over\sinh \kappa\nu
               \sinh \kappa x}
e\sp{ x \kappa\coth\kappa\nu}=
\kappa\left(\coth\kappa x-\coth\kappa\nu\right)e\sp{ x \kappa\coth\kappa\nu}
\nonumber
\\
\label{Phiinf}
\\
\Phi\sp\prime(x;\nu)&=& - \kappa\, \Phi(x;\nu)\left(
\coth\kappa (\nu-x) +\coth\kappa x-\coth\kappa\nu
\right).
\nonumber
\end{eqnarray}
In particular
\begin{equation}
\Phi(x;\infty)={\kappa\over\sinh \kappa x},
\quad{\rm and}\quad
{\Phi\sp\prime(x;\infty)\over \Phi(x;\infty)}=-\kappa\, \coth\kappa x.
\label{Phiinfs}
\end{equation}

\end{document}